\documentclass{amsart}
\usepackage{amsmath,amsthm,amssymb,xypic}
\usepackage[all]{xy}

\theoremstyle{plain}  

\newtheorem{thm}{Theorem}

\newtheorem{lem}[thm]{Lemma}

\newtheorem{rem}[thm]{Remark}
\newtheorem{exa}[thm]{Example}

\theoremstyle{definition}

\newtheorem{dfn}[thm]{Definition}

\theoremstyle{remark}

\newtheorem{claim}{Claim}

\DeclareMathOperator{\Bsl}{Bsl}

\DeclareMathOperator{\bs}{Bsl}

\DeclareMathOperator{\mult}{mult}

\DeclareMathOperator{\Diff}{Diff}
\DeclareMathOperator{\Sing}{Sing}

\DeclareMathOperator{\WCl}{WCl}

\newcommand{\QED}{\ifhmode\unskip\nobreak\fi\quad {\rm Q.E.D.}} 

\newcommand\flip{\dasharrow}
\newcommand\iso{\cong}
\newcommand\f{\varphi}

\newcommand{\F}{\mathbb{F}}  
\newcommand{\C}{\mathbb{C}}
\newcommand{\Z}{\mathbb{Z}}

\newcommand{\I}{\mathcal{I}}

\newcommand{\Q}{\mathbb{Q}}

\renewcommand{\P}{\mathbb{P}}
\renewcommand{\H}{\mathcal{H}} 
\newcommand{\R}{\mathbb{R}}
\renewcommand{\O}{\mathcal{O}}  
\renewcommand{\L}{\mathcal{L}}
\newcommand{\self}[2]{#1\rightleftharpoons#2}
\title[The importance of being $\Q$-factorial]{Birational geometry of quartic 3-folds II:\\
The importance of being $\Q$-factorial}

\author{
Massimiliano Mella}

\address{M. Mella\\
Dipartimento di Matematica\\ 
Universit\`a di Ferrara\\
Via Machiavelli 35\\
44100 Ferrara Italia}
\email{mll@unife.it}
\date{December 2003}
\thanks{Partially supported by EAGER and ``Geometria sulle Variet\`a Algebriche''(MIUR)}
\subjclass{Primary 14E07; Secondary 14J30, 14E30}
\keywords{Mori Fiber spaces, Sarkisov program, birational rigidity, pliability}
\begin{document}
\begin{abstract} The paper explores the birational geometry of terminal quartic 3-folds.
In doing this I develop a new approach to study maximal singularities with positive dimensional
centers. This allows to determine the pliability of a $\Q$-factorial 
quartic with ordinary double points, and it 
shows the importance of $\Q$-factoriality in the context of birational geometry of uniruled 3-folds.
\end{abstract}
\maketitle

\section*{Introduction} 
Let $X$ be a uniruled 3-fold, then $X$ is generically
covered by rational curves. It is a common belief that both
biregular and birational geometry of $X$ are somehow governed by
these families of rational curves. In this paper I am interested in birational
geometry of these objects. The Minimal Model
Program states that such  a $X$ is birational to a Mori fiber
Space (Mfs). Roughly saying after some birational modification either  $X$ can be
 fibered in rational surfaces or rational curves or it becomes
 Fano. For a comprehensive introduction to
this realm of ideas as well as for the basic definitions and results see \cite{CR} and \cite{Co2}. 

In the attempt to tidy up the birational geometry of 3-fold Mori fiber Spaces 
we introduced the notion of pliability, \cite{CMp}.

\begin{dfn}[Corti] If $X$ is an algebraic variety, we define the {\em pliability} of
$X$ to be the set
 \[
\mathcal{P}(X)=\bigl\{ \text{Mfs} \; Y\to T\mid Y\;\text{is birational to} 
\;X \bigr\}/\text{square equivalence}.
 \]
We say that $X$ is {\em birationally rigid} if $\mathcal{P}(X)$ consists of 
one element.
\end{dfn}

It is usually quite  hard to determine the pliability of a given Mori Space, and not many examples
are known. The first rigorous result dates back to Iskovskikh and Manin, \cite{IM}. The main theorem of
\cite{IM} states,
in modern terminology,  that any birational map $\chi:X\flip Y$ from a smooth quartic $X\subset\P^4$
to a Mori fiber space is an isomorphism. This means that $\mathcal{P}(X)=\{X\}$ and
$X$ is birationally rigid.
 
On the other hand consider a quartic 
threefold $X\subset \P^4$ defined by \hbox{$\det M=0$,} where $M$ is a $4\times 4$
matrix of linear forms. One can define a map {\hbox{$f:X\flip\P^3$}} by the assignment
$P\mapsto(x_0:x_1:x_2:x_3)$, where $(x_0,x_1,x_2,x_3)$ is a solution
of the  system of linear equations obtained substituting the coordinates of $P$ in $M$.
For $M$ sufficiently general such a map is well defined and birational. In this case
 $f$ gives a rational parameterization of $X$. The singularities of $X$ correspond to points where
the rank drops.
It is not difficult to show that, for a general $M$, the corresponding quartic has only ordinary double points
corresponding to points where the  rank is 2. Thus a general determinantal quartic threefold 
has only ordinary double points and it is rational. 

From the pliability point of view  this is
discouraging. Minimal Model Theory requires to look at terminal $\Q$-factorial 3-folds and 
ordinary double points are the simplest possible terminal singularities. It would be 
unpleasant if a bunch of ordinary double points were to change a 
rigid structure to a rational variety. 
The point I want to
stress in this paper is that the rationality of a determinantal quartic is due to the lack
of $\Q$-factoriality and not to the presence of singularities. 

\begin{thm} Let $X_4\subset \P_{\C}^4$ be a $\Q$-factorial 
quartic 3-fold with only ordinary double points as singularities. Then
$X$ is neither  birationally equivalent to a conic bundle nor to a
fibration in rational surfaces. Every birational map $\chi:X\flip Y$ to a
Fano 3-fold is a self map, 
that is $Y\iso X$, in particular $X$ is not rational. This is to say
that $X$ is birationally rigid.
\label{th:C}
\end{thm}

\begin{rem}
The case of a general quartic with one ordinary double point has been treated by Pukhlikov, \cite{Pu}. 
Observe that in this case $X$ is automatically $\Q$-factorial. More recently
Grinenko studied the case of a general quartic containing a plane.
\end{rem}

A variety is said to be $\Q$-factorial if every Weil divisor is $\Q$-Cartier. 
Such an innocent definition is quite subtle when realized on a projective variety. 
It does depend both on the kind of
singularities of $X$ and on their position. To my knowledge there are very few papers that tried to
shed some light on this question, \cite{Cl} \cite{We}. In the case of a Fano 3-fold, 
$\Q$-factorial is equivalent to $\dim H^2(X,\Z)=\dim H_4(X,\Z)$, a global topological property, invariant 
for diffeomorphic Fano 3-folds. A recent paper of Ciliberto and Di Gennaro,
\cite{CDG}, deals with hypersurfaces with few nodes. The general behavior is that the presence of
few nodes does not
break $\Q$-factoriality. This is  not true even for slightly worse singularities, as the
following example shows.

\begin{exa}[Koll\'ar] Consider the linear system $\Sigma$, 
of quartics spanned by the following set of monomials
$\{x_0^4,x_1^4, (x_4^2x_3+x_2^3)x_0,x_3^3x_1,x_4^2x_1^2\}$. Then a general quartic
$X\in \Sigma$ has a unique singularity $P$ at $(0:0:0:0:1)$ and the quadratic term is
a general quadric in the linear system spanned by $\{x_3x_0,x_1^2\}$, so that  
analytically $P\in X\sim (0\in (xy+z^2+t^l=0))$ and $P$ is a $cA_1$ point.
The 3-fold  $X$ is not $\Q$-factorial since the plane $\Pi=(x_0=x_1=0)$ is contained in $X$.
The idea is that a general quartic containing a plane has 9 ordinary
 double points, the intersection of the two
residual cubics. In the above case the two cubics intersect just in the point $P$.
\end{exa}

There is a slightly stronger version of Theorem \ref{th:C}.

\begin{thm}\label{th:k} Let $X_4\subset \P_k^4$ be a $\Q$-factorial 
quartic 3-fold with only ordinary double points as singularities over a field $k$, not
necessarily algebraically closed, of characteristic 0. 
Then ${\mathcal P}(X)=\{ X\}$.
\end{thm}

If one considers non algebraically closed fields then peculiar aspects of factoriality and its relation with
birational rigidity appear. Theorem \ref{th:k} and its
significance in this contest, were suggested by J\'anos Koll\'ar.

\begin{exa} Consider the following quartic $Z$
$$(x_0^2+x_1^2)^2+(x_2^2+x_3^2)^2+x_4C=0.$$
Then $Z$ is not $\Q$-factorial over $\C$ because $(x_4=0)_{|Z}$ is a pair of quadrics, say
$Q$ and $\overline{Q}$. 
For a general cubic $C$ the singular points of $Z$ are twelve distinct
ordinary double points. The existence of Minimal Model Program for 3-folds implies that
$Z$ is birational to some Mori Space $Y\neq Z$. Indeed this is the midpoint of a Sarkisov link. 
This can be easily seen with the unprojection method developed by Reid, \cite{un}. The equation of $Z$ is
$$Q\overline{Q}+HC=0.$$
We can introduce the two ratios $y=Q/H=-C/\overline{Q}$ and
$z=\overline{Q}/H=-C/Q$. These are both of degree one and unproject $Z$ to the following complete
intersections
$$\begin{array}{cc}X=\left\{\begin{array}{ll}
 yH= & Q\\ 
 y\overline{Q}=& -C
 \end{array}
 \right.\subset\P^5 &
X'=\left\{\begin{array}{ll}
 zH= & \overline{Q}\\ 
 zQ=& -C
 \end{array}
 \right.\subset\P^5 
\end{array}
$$
 $X$ and $X'$ are projectively equivalent, thus we have a
Sarkisov self link, see \cite{Co2},
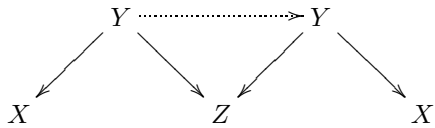
   \[
 \xymatrix{
   &Y\ar[dl]\ar[dr]\ar@{.>}[rr]&&Y\ar[dl]\ar[dr]&  \\
 X&              &Z&&X}
 \]
In the paper I express similar self links with the following compact notation
$$\self{X}{Z_4\subset\P^4}$$ 
In particular the Weil divisors group on $Z$ is generated by
$Q$ and $\overline{Q}$. 
 The two quadrics are conjugated under complex conjugation, so that over $\R$ they are not
defined individually. In particular $Z/\R$ is $\Q$-factorial, hence birationally rigid
by Theorem \ref{th:k}. 
Observe that $X$ is not defined over $\R$.
\end{exa}

Before explaining the proof of Theorem \ref{th:C}, 
let me just give a brief look at the determinantal quartic from the point of view of Sarkisov program.
Let $X=(\det M=0)\subset \P^4$, with $M$ general. Consider a Laplace expansion
of $\det M$ with respect to the $j$-th row. Then the equation of $X$ has the form
$\sum_i l_iA_{ji}=0,$
where the $l_i$ are linear forms and the $A_{ji}$ are cubic forms. Then the
 $A_{ji}$s generate the ideal of a smooth surface $B_j^r$ of degree 6, a Bordiga surface.
It is easy to see that $B_j^r$ passes through all singular points of $X$. The latter are the rank two points
therefore any order three minor has to vanish.
Therefore $X$ is not factorial and consequently not $\Q$-factorial (terminal Gorenstein
$\Q$-factorial singularities are factorial).
The symmetry between rows and columns, in the Laplace expansion,
 suggest that $X$ is a midpoint of a Sarkisov link. Indeed this is the case of a well known 
``determinantal'' involution of $\P^3$, \cite{Pe},
     \[
\self{\P^3}{X\subset\P^4}
 \]

{\em Acknowledgments} This paper found
 his way throughout the darkness of my desktop drawer after
motivating  discussions with Miles Reid during a short visit at Warwick for 
the ``Warwick teach-in on 3-folds'' in January 2002. 
I am deeply indebted with  Alessio Corti for advices, and
much more. His frank criticism on a preliminary version helped to
 improve the paper. The referee, beside numerous suggestions and
 corrections,  found a gap in the first version of
 Lemma \ref{le:ecubo} and suggested a patch. The actual content of
 Lemma \ref{le:nikos} was communicated to me by Nikos Tziolas. 
It is a pleasure to thank  J\'anos Koll\'ar for the  suggestions  and comments he gave me 
during a very pleasant stay in Napoli, for the ``Current Geometry'' Conference 2002.

\section{Maximal singularities and the main theorem}

I start rephrasing Theorem \ref{th:C} in the following form.

\begin{thm}\label{th:main} Let $X_4\subset \P_{\C}^4$ be a $\Q$-factorial 
quartic 3-fold with only ordinary double points as singularities. Then any birational map
$\chi:X\flip V$ to a Mori space $V/T$ is a self map, i.e. $V\iso X$.
\end{thm}

To prove Theorem \ref{th:main} I use the Maximal singularities method combined with Sarkisov
Program, as described in \cite{Co2}, \cite{CPR} and \cite{CM}.
I rely on those papers for the very basic definitions like Mori fiber
spaces, weighted projective spaces, Sarkisov program and links, and philosophical background.
Here I quickly recall what is needed.

 \begin{dfn}[degree of $\chi$] \label{dfn:deg}
 Suppose that $X$ is a Fano 3-fold with the property that $A=-K_X$ generates
the Weil divisor class group: $\WCl X=\Z\cdot A$ (this holds in our case
under the $\Q$-factoriality assumption).
Let $\chi\colon X\dasharrow V$ be a birational map to a
given Mori fiber space $V\to T$, and fix a very ample linear system $\H_V$
on $V$; write $\H=\H_X$ for the birational transform $\chi^{-1}_*(\H_V)$.

The {\em degree} of $\chi$, relative to the given $V$ and $\H_V$, is the
natural number $n=\deg\chi$ defined by $\H=nA$, or equivalently
$K_X+(1/n)\H=0$.
 \end{dfn}

 \begin{dfn}[untwisting] \label{dfn:un}
 Let $\chi\colon X\dasharrow V$ be a birational map as above, and $f\colon
X\dasharrow X'$ a Sarkisov link. We say that $f$ {\em untwists}
$\chi$ if $\chi'=\chi\circ f^{-1} \colon X'\dasharrow V$ has degree smaller 
than $\chi$.
 \end{dfn}

 \begin{dfn}[maximal singularity] \label{dfn:ms}
Let $X$ be a variety and $\H$ a movable linear system.
Suppose that $K_X+(1/n)\H=0$ and $K_X+(1/n)\H$ has not
canonical singularities.  A {\em maximal singularity} is a terminal (extremal) extraction $f:Y\to X$
in the Mori category, see \cite[\S 3]{CM}, having exceptional {\em irreducible}
divisor $E$ such that $f^*(K_X+c\H)=K_Y+c\H_Y$, where $c<1/n$ is the canonical threshold.
The image of $E$ in $X$, or the center $C(X,v_E)$ of the valuation $v_E$, is called the
{\em center} of the maximal singularity.
 \end{dfn}

\begin{rem} In this paper all maximal singularities will be either the blow up 
of an ordinary double point, or generically the blow up of the ideal of a curve $\Gamma\subset X$.
In both cases this is the unique possible maximal singularity with these centers.
This is easy for  curves, while for an ordinary double point it is due to Corti, \cite[Theorem 3.10]{Co2}.
\label{re:uno}
\end{rem} 

\begin{lem}[{\cite[Lemma 4.2]{CPR}}] \label{le:untw} Let $X$, $V/T$, $\H$ be as before,
$\chi: X \dasharrow V$ a birational map. If $E\subset Z \to X$ is a maximal 
singularity, any link $X \dasharrow X'$, starting with the extraction $Z\to X$,
untwists $\chi$.
\end{lem}

The above Lemma, together with Sarkisov program, allow to restrict the attention on
maximal singularities. To study maximal singularities there is an invariant which is very
often useful: the self intersection of the exceptional divisor.
The next Lemma allow to compute $E^3$ when the center is smooth curve
on $X$. To do this I have to determine the correction terms that are
needed to make adjunction formula work in the presence of $cA_1$
singularities along $\Gamma$. This is done using the theory of
Different developed in \cite[\S 16]{U2} and the following Lemma 
kindly suggested by Nikos Tziolas. In the statement and proof of the
Lemma I need a notion of singularity for pairs curve and surface with
$A_t$ points.  
\begin{dfn}
Assume that $p\in S\sim(0\in (xy+z^{t+1}=0))$, for some $t\geq 1$, and
$\Gamma\subset S$ is a smooth curve
through $p$. Let $\nu:U\to S$ be a minimal resolution with exceptional divisors $E_i$,
with $i=1,\ldots,t$. Here I mean that the rational chain starts with $E_1$, ends with $E_t$, and
for $1<i<t$ the intersection $E_i\cdot E_j$ is non zero if and only $j=i\pm1$.

I say that $(\Gamma,S)$ is an $A^k_t$ singularity if $C_U\cdot
E_k=1$ (here and all through the paper I decorate with $_T$ the
strict transform of objects on a variety $T$).
Observe that since $C$ is smooth then $C_U\cdot E_i=0$ for any $i\neq
k$.
\label{df:pair}
\end{dfn}

\begin{lem}[\cite{Tzp}]\label{le:nikos} Let $(0\in X)$ be a $cA_1$ singularity and
  $0\in\Gamma\subset X$ a smooth curve through it. Let $f:Y\to X$ be a
  terminal  extraction with center a 
smooth curve $\Gamma$
and exceptional divisor $E$. Then $f$ can be obtained from the diagram
$$
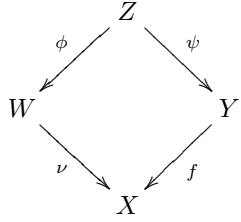
\xymatrix{
&Z\ar[dl]_{\phi}\ar[dr]^{\psi}&\\
 W\ar[dr]_{\nu}&&Y\ar[dl]^{f}\\
&X&}
$$
such that
\begin{itemize}
\item[i)] $W$ is the blow up of $X$ along $\Gamma$. The
  $\nu$-exceptional divisors are a ruled surface $E$ over $\Gamma$ and
  $F\iso\P^2$ aver the singular point. $Z$ is a $\Q$-factorialization
  of $E$ and $\psi$ contracts $F_Z\iso F$ to a point.
\item[ii)] $S_Y=f^{-1}_*S\iso S$, where $S$ is a general  section of
  $X$ through $\Gamma$.
\item[iii)] $(\Gamma,S)$ is an $A^k_{2k-1}$ singularity.
\end{itemize}
\end{lem}
\begin{proof} First prove that $W$ is $cA$, \cite{Ko}.
Let $S$ be the general section of $X$ through $\Gamma$. Then one can
assume, \cite{Tza}, that $S$ is given by $xy-z^{n+m}=0$ and $\Gamma$
by $x-z^n=y-z^m=0$, for some $n\leq m$, equivalently $S$ by
$xy+xz^n+yz^m=0$ and $\Gamma$ by $x=y=0$. Then $X$ has the form
$$xy+xz^n+yz^m+tg_1(x,y,z,t)+tg_{\geq 2}(x,y,z,t)=0$$
and $\Gamma$ is $x=y=t=0$. To have a $cA_1$ singularity the quadratic
term
$ xy+tg_1(x,y,z,t)$ must be irreducible. Now a straightforward
explicit computation of the blow up of the maximal ideal of $\Gamma$
shows that $W$ is $cA$.

Then by \cite{Tzb} it follows that $Z$ and hence $Y$, can be constructed in
families. Therefore we may study the deformed equation 
$$xy+xz^n+yz^m+tg_1(x,y,z,t)+t^k=0$$
for $k\gg1$. The blow up computation and the irreducibility of the
quadratic term yields that $W$ has isolated singularities along $E\cap
F$. Therefore $Z$ is just the blow up of $E$ and hence $F_Z\iso
F\iso\P^2$. This also proves that $F_Z$ is contracted to a point by
$\psi$. 

To see the claim on $S_Y$ take a general member $S_W\in|-K_W|$. Then
$S_W$ has $A_i$ singularities and avoids the singular points along
$E\cap F$. Let $C=S_W\cap F$. The $C$ is contracted by $\nu$ and
therefore $S_Y=\nu(S_W)\iso S$.

Since $W$ is smooth on the generic point of $E\cap F$, \cite[Proposition 4.6]{Tza}, it follows that
$(\Gamma, S)$ is an $A^k_{2k-1}$ singularity
because in any other case $W$ would be singular at $E\cap F$.
\end{proof}

Next I derive the numerical result about self intersection
from Lemma \ref{le:nikos}.

\begin{lem} Let $f:Y\to X$ be a terminal extraction with center a 
smooth curve $\Gamma$
and exceptional divisor $E$. Assume that $X$ has only $cA_1$ points
along $\Gamma$. 
Let $\Sigma$ be any linear system with $\Bsl\Sigma=\I_{\Gamma}$ and $S\in \Sigma$ 
a general element. Assume that $S$ is normal. Then $f_{|S_Y}:S_Y\to S$ is an isomorphism,
 and
$$E^3=-S\cdot\Gamma-(\Gamma\cdot\Gamma)_S$$
or equivalently
$$E^3=K_X\cdot\Gamma-2g(\Gamma)+2-\Diff(\Gamma,S)$$
\label{le:ecubo}
\end{lem}

\begin{rem}
In the hypothesis of Lemma \ref{le:ecubo} one can define 
the different of $\Gamma$ in $X$ as
$$\Diff(\Gamma,X):=K_X\cdot\Gamma-E^3-2g(\Gamma)+2$$
This suggests the possibility to extend the theory of Different, \cite[\S 16]{U2}, 
to higher codimension subvarieties.
\end{rem}

\begin{proof}We already proved, Lemma \ref{le:nikos}, that $S_Y\iso
  S$. 
 By hypothesis $f^*(S)=S_Y+E$, and consequently
$$E^3=f^*(S)\cdot E^2-S_Y\cdot E^2$$
Projection formula yields $f^*S\cdot E^2=-S\cdot \Gamma$. By Lemma
\ref{le:nikos}   $E_{|S_Y}=\Gamma$, therefore $S_Y\cdot E^2=
(E_{|S_Y}\cdot E_{|S_Y})_{S_Y}=(\Gamma\cdot\Gamma)_S$.
Note that  $K_S=(K_X+S)_{|S}$, therefore
$$(\Gamma\cdot \Gamma)_S=2g(\Gamma)-2-K_X\cdot\Gamma-S\cdot\Gamma+\Diff(\Gamma,S)$$
by
adjunction formula compensated by the Different, \cite[Chapter 16]{U2}. 
\end{proof}

I now go back to Theorem \ref{th:main}. 
The first task is to recognize birational maps.
The geometry of $X$ suggests the existence of some birational self maps, the ``Italian'' approach,
according to \cite{CPR}: 
\begin{itemize}
\item the reflection through a singular point $p$
\item the elliptic involution associated to a line $l$ 
containing some singular point.
\end{itemize} 
The general line through $p$ intersect the quartic
in two more points $Q_1$ and $Q_2$. The self map suggested is $Q_1\mapsto Q_2$. 
A general plane containing $l$ has a smooth cubic $C$ as residual intersection with $X$.
Furthermore a singularity, say $P$,
 provides the family of these cubics of a section, namely a common origin
to the group structure. The self map suggested is $R\mapsto -R$ where $-R$ is the inverse of
$R$ in the group structure on $C$ with origin $P$.

Then I describe those maps in terms of Sarkisov links.

\noindent After \cite{Co2}, \cite{CPR} and \cite{CM} this is now a nice and pleasant exercise.
Indeed the only possibility that is not yet described in neither \cite{Co2} nor \cite{CM} is the
one of a line with three singularities along it.
Assume that $l\subset X$ is a line with three distinct singular points along it. Note that
this is the maximum number of singular points along a line on a
quartic  with isolated singularities.
After a coordinate change we can assume that $l=(x_2=x_3=x_4=0)$ and 
the equation of $X$ has the following form
$$F=L(x_0^2x_1+x_1^2x_0)+Q_1x_0^2+Q_2x_1^2+Q_3x_0x_1+C_1x_0+C_2x_1+D$$

Let $f:Y\to X$ be the unique terminal extraction with center $l$ and
exceptional divisor $E$.
I want to understand the anticanonical ring of $Y$. Let
$$H_i=(x_i=0)_{|X}$$
It is immediate that $H_{iY}\in|-K_Y|$. Since $l=(x_2=x_3=x_4=0)$ and
$f$ is generically the blow up of the maximal ideal then $-K_Y$ is
nef. The linear form $\tilde{H}:=L_{|X}$ has multiplicity two along
$l$. Therefore a general plane section of $\tilde{H}$ through $l$ has
residual intersection a conic, say $C$, that, generically, intersects $l$ in two
points. In particular $C_Y\cdot K_Y=0$ and  $NE(Y)=\langle
e,C\rangle$, where $e\subset E$ is $f$-exceptional. 
Note that the special hyperplane section
$\tilde{H}_Y$ is covered by curves proportional to $C$.
Therefore the ray spanned by $[C]$ is not small and by the two ray
game I conclude that there  is no Sarkisov link starting from
the extraction $f:Y\to X$. This is usually called a bad link,
\cite{CPR}, \cite{Co2}.

The only Sarkisov links that have center either a singular point or a
line through a singular point are therefore the following:
\begin{itemize}
\item[$\rho_x$] for any singular point $x\in X$ 
   \[
\self{X}{Z_6\subset\P(1,1,1,1,3)}
 \] 
\item[$\f^l_1$] for any line $l\subset X$ passing through one singular point 
    \[
\self{X}{Z_{12}\subset\P(1,1,1,4,6)}
 \] 
\item[$\f^l_2$] for any line $l\subset X$ passing through two singular points 
    \[
\self{X}{Z_{8}\subset\P(1,1,1,2,4)}
 \]
\end{itemize}

Note that to a line with more than one singularity are associated different elliptic involutions.
I can choose any singular point as origin on the elliptic curves. 
But still, by Sarkisov theory, the maximal singularity with center the  line is unique.
This is because the elliptic involution, in this case, is a composition elementary links.

To prove Theorem \ref{th:main} it is now enough 
 to show that any  birational map can be factored by the self maps
 described. 
It is now standard, see  \cite[\S 3]{CPR}, that this is equivalent
to the following.
\begin{thm} Let $X_4\subset\P_{\C}^4$ be a $\Q$-factorial quartic 3-fold with only ordinary double points and
$E$ a maximal singularity. Then either:
\begin{itemize}
\item[-] the center $C(X,v_E)=p$ is a singular point, or
\item[-] the center $C(X,v_E)=l$ is a line through some singular point.
\end{itemize}
In both cases the assignment
 identifies the maximal singularity, hence the Sarkisov link, uniquely, see Remark~\ref{re:uno}.
\label{th:center}
\end{thm}

The proof of Theorem \ref{th:center} is the core of the next section.

\section{Exclusion}

A maximal center on a Fano 3-fold is either a point or a curve. 

\noindent The case of smooth points can be treated with many different techniques.
The main result of \cite{IM} is indeed that  a smooth point is not a maximal center on a quartic. 
Corti, \cite{Co2}, gave an amazingly simple proof
using numerical properties of linear system on surfaces. The recent classification of Kawakita,
\cite{Kw}, gives a third possible proof based on terminal extractions, see \cite[Conjecture 4.7]{Co1}.

I am therefore bound to study centers of positive dimension. I can actually prove a stronger statement.

\begin{thm}\label{th:cA1} Let $X_4\subset \P_{\C}^4$ be a $\Q$-factorial 
quartic 3-fold with only $cA_1$ points. Assume that a curve $\Gamma$ is a center of 
maximal singularities. Then $\Gamma$ is a line through some singular point.
\end{thm}
 
\begin{rem} The Theorem is an important step in the direction 
 of \cite[Conjecture 1.3]{CMp} 
\end{rem}

\begin{proof}
From now on $\Gamma\subset X$ will be an irreducible curve assumed to be the  center of a 
maximal singularity.
The unique terminal extraction is then generically the blow up of the ideal of $\Gamma$ in $X$.
Therefore the linear system $\H\subset|\O(n)|$, associated to the extraction, satisfies
 $$\gamma=\mult_{\Gamma}\frac1n\H> 1.$$ 
We prove the theorem in several steps
\begin{description}
\item[Step 1] A raw argument shows that $\deg \Gamma \leq 3$.
\item[Step 2] $\Gamma$ can not be a space curve.
\item[Step 3] If $\Gamma$ is a plane curve then it is a line through some singular point.  
\end{description}

\paragraph*{\textsc{Step 1}} Choosing general members $H_1$, $H_2$ of $\H$ and 
intersecting with a general hyperplane section $S$ we obtain
 \[
 4n^2=H_1\cdot H_2\cdot S>\gamma n^2\deg\Gamma.
 \]
This implies that $\deg \Gamma \leq 3$.

\paragraph*{\textsc{Step 2: space curves}} If $\Gamma$ is a space curve, then
by Step~1 it must be a rational normal curve of degree 3, contained in
a hyperplane $\Pi \cong \P^3 \subset \P^4$.
Let $S \in |I_{\Gamma,X}(2)|$ be a general quadric 
vanishing on $\Gamma$,$\L$ the mobile part of $\H_{|S}$; write 
 \[
 \O_S(1)=\frac{1}{n}\H_{|S}=L+\gamma\Gamma,
 \]
where $L=(1/n)\L$ is nef. Note that, because $I_\Gamma$ is cut out by
quadrics, 
\begin{displaymath}
  \mult_\Gamma \H = \mult_\Gamma \H_{|S} = n \gamma > n.
\end{displaymath}

Let $f:Y\to X$ be the maximal singularity, with exceptional divisor $E$, 
and center $\Gamma$. By Lemma \ref{le:ecubo} we can compute $E^3$ by means of
$\Diff(\Gamma,S)$.
Assume that $(C,U)$ is an $A^k_t$ singularity, keep in mind Definition
 \ref{df:pair}, let
 $\nu^*(C)=C_W+\sum d_i E_i$ then it is a straightforward check 
on the intersection matrix of an $A_t$
singularity,
see for instance \cite[pg 16]{Ja}, that
\begin{equation}
\label{eq:akn}
(t+1)d_i=\left\{\begin{array}{ll}
 i(t-k+1) & {\rm\   if\ } i\leq k \\ 
 (t-i+1)k& {\rm\   if \ }i\geq k
 \end{array}
 \right.
\end{equation}
Incidentally observe that $\Diff(C,U)_x=C_W\cdot\sum d_iE_i=d_k$.

I now come back to our original situation by Lemma \ref{le:nikos} part
iii) $(\Gamma,S)$ is an  $A^k_{2k-1}$ singularity for
some $k$, with $l/2\geq k\geq 1$. In particular the Different is
\begin{equation}
\label{eq:k2}
\Diff(\Gamma,S)_p=k/2
\end{equation}

This proves, together with Lemma \ref{le:ecubo}, that 
\begin{equation}
\label{eq:ecubo}
E^3=-3+2-\sum_{p_i\in \Sing(S)} k_i/2
\end{equation}

Then we need to bound the contributions of the singularities globally.

\begin{lem}In the above notation  $\sum_{p_i\in \Sing(S)} k_i\leq 7$.
\label{le:bound}
\end{lem}

\begin{proof}
To prove the bound we need the following reinterpretation of the $k_i$'s, see also
\cite{Tza}. By Lemma \ref{le:nikos} part iii)
we can realize $p_i\in\Gamma\subset S\subset \Q^3$ analytically as 
$$0\in(x=y=0)\subset (xy+yz^{k_i}+xz^{k_i}=0)\subset \C^3,$$ 
see for instance \cite[pg 13]{Ja}. Let $\mu_i:Z\to\C^3$ be the 
 blow up of $(x=y=0)$, with exceptional divisor
$E_Z$ and $F_i=\mu_i^{-1}(0)$. Then $S_{Z|E_Z}=k_iF_i+{\rm effective}$. 
Let $\nu:W\to \P^4$ be the blow up of $\Gamma$, with exceptional divisor $E_W$,
and $F_i=\nu^{-1}(p_i)$. Since $\Bsl|\I_{\Gamma,\P^4}(2)|=\Gamma$  then 
$X_{W|E_W}=k_i F_i+{\rm effective}$. 
For any divisor $D\subset\P^4$ such that $D\supset\Gamma$  and $D_{|X}$ is smooth on the generic point 
of $\Gamma$ we have
$$(D_Z\cap X_Z)_{|E_Z}=h_i F_{i|D_Z}+{\rm effective},$$
 for some 
\begin{equation}
\label{eq:quad}
h_i\geq k_i
\end{equation}
Thus to bound the global
contribution it is enough to understand the normal bundle of $\Gamma$ in some smooth
divisor $D$ such that $D_{|X}$ is smooth on the generic point of $\Gamma$.
Let $H=\Pi_{|X}$ be the unique hyperplane section containing $\Gamma$, and
$H_{|S}=\Gamma+\Delta$.
I claim that $\Delta\not\supset\Gamma$. 
Assume the opposite and let $H_{|S}=2\Gamma+R$. Then $\deg R=2$ and
$R$ is a pair of skew lines, say, $l_1$, $l_2$,
secant to $\Gamma$. Since $S$ is general then  $l_i\not\subset \Bsl\H$ and 
$l_i\cap\Gamma\not\subset (\Sing(X)\cap\Gamma)$. 
Then we derive the impossible
$$1=\H/n\cdot l_i\geq\gamma\Gamma\cdot l_i>\gamma.$$
We can therefore choose $D=\Pi$.
It is well known, \cite{Hu}, that  $N_{\Gamma/\P^3}\iso\O(5)\oplus\O(5)$.
Let $\nu:W\to D$ be the blow up of $\Gamma$ with exceptional divisor $E_W\iso \F_0$.
Then $X_{W|E_W}\equiv f_0+7f_1$, where $f_1$ is a fiber of $\nu$.
The inequality $\sum k_i\leq \sum h_i\leq7$ is obtained.
\end{proof}

Consider again, the hyperplane section  $H=\Pi_{|X}$ and the maximal singularity $f:Y\to X$. 
Let $D_Y\subset Y$ be any effective irreducible divisor distinct from $E$. 
Then $D_Y=f^*D-\alpha E$, for some positive $\alpha\in\Q$ and $D\in |\O(d)|$. 
The divisor $D_Y$ is numerically equivalent to $d H_Y+(d-\alpha)E$.  To conclude the step it is,
therefore enough
to prove that the cone of effective divisors on $Y$ is generated by $H_Y$ and $E$.
This is the content of the next Lemma.

\begin{lem} $NE^1(Y)=\langle H_Y,E\rangle$.
\label{le:cono}
\end{lem}

\begin{rem} This is just a rewriting of the usual exclusion trick.
I prove that a linear system like $\H$ has to have a fixed component,
in this case $H$. I hope that in this way  it is easier to digest and maybe 
generalize. See also Remark \ref{rem:prob}.
\end{rem}

\begin{proof}
Let $B_Y\subset Y$ be any effective irreducible $\Q$-divisor distinct from $E$ and $H_Y$. 
Then $B_Y=f^*B-\beta E$, for some positive $\beta\in\Q$ and $B\in
|\O(b)|$. Actually $\beta\in \Z$  since $X$ has index 1 and is $\Q$-factorial. I have 
 to prove that $\beta\leq b$. By Lemma \ref{le:ecubo} $\dim \Bsl|S_Y|\leq 0$, hence
the cycle $B_Y\cdot H_Y\cdot S_Y$ is effective.
The following inequality is satisfied
$$\begin{array}{rl}
0\leq& 
B_Y\cdot H_Y\cdot S_Y=(f^*(B)-\beta E)(f^*(H)-E)(f^*S-E)\\
=& 8b-3b-9\beta-\beta E^3=5b-(8-\sum k_i/2)\beta
\end{array}
$$
This proves the claim for $\sum k_i\leq 6$.

Assume that $\sum k_i=7$. First I need to better 
understand this special configuration of singularities.

Let $\nu:W\to \Pi$ be the blow up of $\Gamma$ with exceptional divisor $E_W\iso \F_0$,
 $g$ an ``horizontal'' ruling of $E_W$, and $f_i$ fibers of $\nu$. Then the assumption on singularities 
yield
$H_{W|E_W}=g+\sum_1^7 f_i$, where the $f_i$ are not necessarily
 distinct. Note that for each point $y\in\Gamma$ there is a quadric
 cone $Q^y\subset\Pi$ containing $\Gamma$ and with vertex $y$. Then
 $Q^y_{W|E_W}=g_y+f$. In particular for any $g\subset E_W$ ``horizontal''
 ruling there exists a
 quadric cone $Q^g\subset \Pi$ such that $Q^g_{W}\supset g$.
This proves that there exists a quadric cone, say $\tilde{Q}\subset\P^3$,
such that $\tilde{Q}_{|H}=2\Gamma+C$, for some conic $C$. Similarly there exists
a cubic surface $\tilde{M}$ such that $\tilde{M}_{|H}=2\Gamma+R$ and $\tilde{M}_{|\tilde{Q}}=2\Gamma$.
Therefore the equation of $H$ can be written as
$$\tilde{Q}K+\tilde{M}P=0,$$
where $K$ is a quadric and $P$ is a linear form. 

Assume that $\Pi=(x_4=0)$. Let $\Sigma$ be the linear system of quadrics spanned by 
$\{\tilde{Q},x_4x_0,\ldots,x_4x_3\}$. Fix $\overline{S}\in\Sigma_{|X}$ a general element.
By construction we have $H_{|\overline{S}}=2\Gamma+C$. Let 
$$H_{Y|E}=\Gamma_0+F,$$
where $F$ is $f$-exceptional.
Then for effective,
$f$-exceptional divisors $F'$ and  $G$, we have 
$$\overline{S}_{Y|E}=\Gamma_0+F^{\prime},\ \ \overline{S}_{Y|H_Y}=\Gamma_0+C+G$$
and
\begin{equation}
\label{eq:rel}
(\overline{S}_Y-E)\cdot H_Y=C+G-F
\end{equation}
\begin{claim} $F-G\equiv\O_E$\label{cl:caga}
\end{claim}
\begin{proof}[Proof of the Claim] The cycle $F$ is the $f$-exceptional
 part of $H_{Y}$. The cycle $G$ is $f$-exceptional and it is 
contained in $\Bsl\Sigma_{Y}$ therefore $F-G$ is effective.
Let $\phi:=f_{|E}:E\to\Gamma$ be the restriction morphism and $E^0=E\setminus\phi^{-1}(\Sing(X)\cap\Gamma)$.
In our notation we have $\Bsl(\Sigma_{Y|E})=\Gamma_0+G$, thus we can assume that
$F'= G+M$, for some divisor $M=\phi^*A$ supported on $E^0$. 
 Let me interpret this divisor in a different way. Let $\overline{Q}\in\Sigma$ be the quadric whose
 to $X$ is $\overline{S}$. Since $\tilde{Q}=\overline{Q}_{|\Pi}$ is a cone then
$N_{\Gamma/\overline{Q}}\iso\O(2)\oplus\O(5)$. 
Let $\nu:W\to\overline{Q}$ be the blow up of $\Gamma$, with exceptional divisor $E_W$. Then a computation
similar to that of  Lemma \ref{le:bound} yields
$$X_{W|E_W}\equiv \Gamma_0+\nu_{|E_W}^*\O(10)\ {\rm and\ }X_{W|E_W}=\Gamma_0+\sum h_i f_i+{\rm effective} $$
where the $h_i\geq k_i$ and $\nu(f_i)\in\Sing(X)$.  This proves that
$\deg A\leq 10-\sum k_i=3$. Taking into account a reducible
quadric in $\Sigma$, we have 
$F'\equiv F+\phi^*\O(3)$. This shows that $F-G\equiv \phi^*(A-\O(3))$  and together with the bound
on the degree of $A$ the desired $F-G\equiv\O_E$.
\end{proof}
Projection formula and equation (\ref{eq:ecubo}) at page \pageref{eq:ecubo} yield
$$(\overline{S}_Y-E)\cdot H_Y\cdot E=6-2H_Y\cdot E^2=6-2(f^*H\cdot E^2-E^3)=3$$
Then by Claim \ref{cl:caga} and equation (\ref{eq:rel})  I derive
\begin{equation}
\label{eq:3b}
E\cdot C=3
\end{equation}
Note that if  $C$ is reducible, then each irreducible component $C_i$ is a line. In this case 
the inequality
 $E\cdot C_i\leq 2$ is immediate. Thus we proved that for any irreducible component $C_i\subset C$
\begin{equation}
\label{eq:finalmente}
E\cdot C_i\geq \deg C_i
\end{equation}

Let us assume that $C$ is irreducible, the reducible case is similar and left to the reader.

Assume that $B_{Y|H_Y}=aC+\Delta$, for some effective divisor $\Delta$, with $\Delta\not\supset C$.
The above construction gives
$$(f^*(B-a\overline{S})-(\beta-2a)E)_{|H_Y}=\Delta+a(F-G)$$
This proves that
$(f^*(B-a\overline{S})-(\beta-2a)E)\cdot C\geq 0$
and we conclude  by equation (\ref{eq:finalmente}) that
$$(b-2a)\geq (\beta-2a)$$
\end{proof}

\paragraph*{\textsc{Step 3: plane curves}} Here we assume that
$\Gamma$ is a plane
curve of degree $d$ (by Step~1, $d\leq 3$), other than a line passing 
through some singular point. Let $\Pi \subset \P^4$ be the plane spanned by $\Gamma$.
Fix $S$,
$S'$ be general members of the linear system $|\I_{\Gamma,X}(1)|$. Here it is helpful and convenient to 
treat two cases, namely:
\begin{description}
\item[Case 3.1] $\Gamma\cap \Sing(X)=\emptyset$, and $1\leq d \leq 3$.
\item[Case 3.2] $\Gamma\cap \Sing(X)\neq\emptyset$ and $2\leq d \leq 3$.
\end{description}

\subparagraph*{\textsc{Case 3.1}} I first deal with the easy, and well
known, case of curves in the smooth locus.  
Let $f:Y\to X$ be the maximal singularity  with center $\Gamma$, and
exceptional divisor $E$. Then $Y$ is just the blow up of $\I_{\Gamma}$
and by Cutkosky's classification, \cite{Cu},
of terminal extraction   
$$E^3= K_X\cdot\Gamma-2p_a(\Gamma)+2$$
 
\begin{lem} $NE^1(Y)=\langle S_Y,E\rangle$
\end{lem}
\begin{proof} 
Let $B_Y\subset Y$ be any effective irreducible $\Q$-divisor distinct from $E$ and $S_Y$. 
Then $B_Y=f^*B-\beta E$, for some positive $\beta\in\Z$ and $B\in |\O(b)|$. The claim is
equivalent to prove that $\beta\leq b$.
Consider a general element $ D\in
|I_{\Gamma, X}(d)|$. The cycle $B_Y\cdot S_Y\cdot D_Y$ is effective, thus
$$\begin{array}{rl}
0\leq& 
B_Y\cdot S_Y\cdot D_Y=(f^*B-\beta E)(f^*S-E)(f^*B-E)\\
=&4bd-bd-d\beta -d^2\beta-\beta E^3=3bd-d^2\beta+(2p_a(\Gamma)-2)\beta
\end{array}
$$
It is a simple check that for any possible pair $(d,p_a(\Gamma))$ the equation gives
$\beta\leq b$.
\end{proof}
The Lemma finish off the Case~3.1.

\subparagraph*{\textsc{Case 3.2}} From now on we assume that 
there are singular points along $\Gamma$ and $\Gamma$ is
not a line.

We work with the linear system $\Sigma=|S, S'|$, even though
$\Gamma$ is usually only a component of its base locus $C=S\cap S' =
\bs \Sigma= X \cap \Pi$. Write
 \[
 C = \mu \Gamma + \sum \mu_i \Gamma_i
 \]
We are assuming that $X$ is $\mathbb{Q}$-factorial. This implies that $\Pi$ can not
be contained in $X$, and $C$ is a \emph{curve}. Assume first that 
the intersection $S \cdot S'$ is reduced then 
$\mult_\Gamma \H =\mult_\Gamma \H_{|S}$ and $\mult_{\Gamma_i} \H =\mult_{\Gamma_i} \H_{|S}$.
\label{pag:plane}
We always restrict to $S$ and write
$$
\begin{array}{rl}
A :=(1/n)\H_{|S}&= L+\gamma\Gamma+\sum \gamma_i\Gamma_i \\
       S_{|S}^{\prime}&= C=\Gamma+\sum \Gamma_i
\end{array}
$$
The technique consists in selecting a ``most favorable'' component of
$C$, performing an intersection theory calculation using that $L$
is nef, and get that $\gamma\leq 1$. This inequality
contradicts the hypothesis that  $\Gamma$ is a maximal
singularities. Indeed,
 keeping 
in mind Remark \ref{re:uno}, we have
$$n<\mult_{\Gamma}\H=\mult_{\Gamma}\H_{|S}=\gamma n\leq n$$

\begin{rem}
\label{rem:prob}
This is similar to what I did with the twisted cubic
with $\sum k_i=7$, see Lemma \ref{le:cono}. I believe that the $E^3$ approach works also
in this case, but I did not check it. On the other hand each different configuration needs
different calculations. For this reason I developed a unified approach with more
emphasis on the intersection theory on $S$.
\end{rem}

Because $\Gamma$ is a center of a \emph{maximal}
singularity, $\gamma \geq \gamma_1$, $\gamma_2$, hence possibly after
relabeling components of $C$, we can assume that: 
 \[
 \gamma\geq \gamma_2\geq \gamma_1.
 \]
Consider now the effective $\Q$-divisor
 \[
 (A-\gamma_1 S')_{|S}=L+(\gamma-\gamma_1)
 \Gamma+(\gamma_2-\gamma_1)\Gamma_2.
 \]

I now show that $(\Gamma\cdot\Gamma_1)_{S}\geq \deg \Gamma_1$; together with
the last displayed equation this implies that $\gamma \leq 1$ and
finishes the proof.

The curve $\Gamma_1$ is either a line or a conic. Let $D\in|\I_{\Gamma_1,\P^4}(\deg
\Gamma_1)|$ a general
element, and $D_{|S}=\Gamma_1+F$. Since $D\cap\Pi=\Gamma_1$ is a complete intersection, then
 $F$ intersects  $\Gamma$ only at
$\Sing(X)\cap\Gamma_1\cap\Gamma$.
Fix a point $p\in \Sing(X)\cap\Gamma\cap\Gamma_1$. Assume that $p\in X\sim (0\in (xy+z^2+t^l=0))$.

Let $f_1:X_1\to X$ be the blow up of $p\in X$ with exceptional divisor $E_1$.
Then $S_{1|E_1}=C_1$ is a conic and since $\Bsl\Sigma=\Pi$ then $C_1$ is reduced.
This proves that  either $S_{1|E_1}$ is smooth or it has one singular point only, say $x_1$, and
$C$ is a pair of lines. Let $f_2:X_2\to X_1$ be the
blow up of $x_1$, with exceptional divisor $E_2$. 
If $p_1\in X_1$ is a smooth point then $E_2$ is a plane, and $\Bsl\Sigma_2$ is contained in a line.
The surface $S_2$ is smooth and 
already $S_1$ was non singular. Otherwise $p_1\in X_1\sim (0\in (xy+z^2+t^{l-2}=0))$ and
we simply repeat the same argument.
This gives a morphism $\nu:W\to X$, with exceptional divisors $G_i$, for $i=1,\ldots,g$.
Such that $\nu_{|S_W}:S_W\to S$ is a minimal resolution. Moreover
$S_W\cap G _i=L_i\cup R_i$ is a pair of disjoint (-2)-curves, for any $i<g$, and $S_W\cap G_g=T$ is
either a (-2)-curve or a pair of (-2)curves intersecting in a point.
This proves that $p\in S_1$ is an $A_m$ point, with $m\leq l$. 
Furthermore $F$ is smooth at $x$.

Number all irreducible components of the  resolution $\nu_{|S_W}$ 
from 1 to $m=2g-\epsilon$, where $\epsilon=1,0$,
according to the parity of $m$. Start with $L_1=:E_1$, then $L_i=:E_i$ and $R_i=E_{m+1-i}$ for any $i<g$.
Similarly let $T=L_g\cup R_g=:E_g\cup E_{g+1}$, where $E_g\cdot E_{g+1}=1$,
 if it is reducible and $T=E_g$ 
if it is irreducible.

As our aim is to calculate an intersection product we need to understand
 the pairs $(\Gamma,S)$, $(\Gamma_1,S)$, and $(F,S)$.

If $(\Gamma_1)_W\cap T=\emptyset$ then there exists an index $j<g$ such that
$(\Gamma_1)_W\cdot E_j=1$ and $F_W\cdot E_{m+1-j}=1$.
If $(\Gamma_1)_W\cap T\neq\emptyset$ and $T=L_g\cup R_g$ is reducible we labeled the component
in such a way that $(\Gamma_1)_W\cdot E_g=1$ and $F_W\cdot E_{g+1}=1$.
Finally for $(\Gamma_1)_W\cap T\neq\emptyset$ and $T$ is irreducible then
$(\Gamma_1)_W\cdot E_g=F_W\cdot E_g=1$. 

In any case $(\Gamma_1,S)$ is of type $A_m^j$, for some $j\leq m+1-j$. While
$(F_Z,S)$ is of type  $A_m^{m+1-j}$.

Let  
$$\nu_{|S_W}^*(\Gamma_1)=(\Gamma_1)_W+\sum r_i E_i$$ 
and 
$$\nu_{|S_W}^*(F)=F_W+\sum f_i E_i.$$

Then the $r_i$s and the $f_i$s are completely determined by equation (\ref{eq:akn}) at page \pageref{eq:akn}.
The index $j$ satisfies the inequality  $j\leq m+1-j$ by hypothesis. Assume that $i\leq m+1-i$ is
also true, then $m+1-j\geq i$. Thus for any index $i$ such that $i\leq m+1-i$
we have,
$$
(m+1)(r_i-f_i)=\left\{\begin{array}{ll}
 (m+1-2j)i & {\rm\   if \ }i\leq j \\ 
 (m+1-2i)j& {\rm\   if \ }i\geq j
 \end{array}
 \right.
$$
These yield
\begin{equation}
\label{eq:conto}r_i\geq f_i \mbox{\rm \ for any $i\leq m+1-i$.}
\end{equation}
The curve $\Gamma\subset\Pi$ has at most a simple node or a simple cusp then 
\begin{equation}
\label{eq:left}
\Gamma\cdot E_i=0 \ \mbox{for any $i>g$ i.e. $i>m+1-i$}
\end{equation}

By construction $F_W\cdot \Gamma_W=0$ therefore by projection formula
 $$(\Gamma_1\cdot\Gamma-F\cdot \Gamma)_x\geq \sum_i(r_i-f_i)\Gamma_W\cdot E_i.$$
By equation (\ref{eq:left}) we can restrict the summation on indexes
satisfying $i\leq m+1-i$ and equation (\ref{eq:conto}) yields
$$(\Gamma_1\cdot\Gamma)_x-(F\cdot \Gamma)_x\geq 0.$$
Finally all contributions coming from singular points give
$$\deg \Gamma_1 \deg\Gamma=D\cdot \Gamma=((\Gamma_1+F)\cdot\Gamma)_{S}\leq 2
(\Gamma_1\cdot\Gamma)_{S},$$
and consequently  the needed bound since $\deg\Gamma\geq 2$.

Next we consider the case in which $S_{|S}^{\prime}=\Gamma+2l$, where $\Gamma$ is a conic and $l$ is a line.
Again $S$ is smooth at $(\Gamma\cap\Gamma_1)\setminus(\Sing(X)\cap\Gamma)$ as well as on the
generic point of $l$.
Indeed we are just fixing a plane,
therefore we can always choose an hyperplane containing $\Pi$, and not tangent to $X$ at both
$(\Gamma\cap\Gamma_1)\setminus(\Sing(X)\cap\Gamma)$ and at the generic point of $l$.
Then $\H_{|S}=\L+\gamma\Gamma+\alpha l$. Consider the $\Q$-divisor
$$(\H-(\alpha/2)S')_{|S}=\L+(\gamma-\alpha/2)\Gamma.$$
Then 
$$(1-(\alpha/2))\geq (\gamma-\alpha/2)\Gamma\cdot l.$$
To exclude this case we argue exactly as before that $\Gamma\cdot l\geq 1$. Keep in mind that also in 
this case 
$S$ has only isolated singularities. Therefore locally around $x$ all the calculations are the same.

Finally we have to treat the double conic case. That is assume that $S_{|S}^{\prime}=2\Gamma$.
If there exists an hyperplane section $\tilde{S}$ such that $\mult_{\Gamma} \tilde{S}=2$ then for a general
hyperplane section $H$
$$4=H\cdot\frac{\H}n\cdot\tilde{S}\geq \frac4n\mult_{\Gamma}\H=4\gamma.$$
We can therefore assume that the tangent space to $X$ along $\Gamma\setminus(\Gamma\cap \Sing(X))$ is not
fixed. It is immediate to observe that for any smooth point $p\in\Gamma$ the embedded tangent space
contains $\Pi$.
Let us assume the following notations:
\begin{itemize}
\item[-] $\Gamma\subset\Pi\subset\P^4\sim (x_0x_4+x_3^2=0)\subset (x_1=x_2=0)\subset\P^4$,
\item[-] $x\equiv(1:0:0:0:0)\not\in \Sing(X)$,
\item[-] $T_xX=(x_1=0)$,
\item[-] $S=(x_1=0)_{|X}$,
\item[-] $S'=(x_2=0)_{|X}$.
\end{itemize}

By construction $\H_{|S}=g\Gamma+\L$, where $\L$ is a linear system without fixed components and
$g\geq \gamma$.
Up to consider $2\H$ we can further assume that
$g=2k$ is even.  Since $S\cdot S'=2\Gamma$ then
 a general divisor $H\in\H$  has an equation of type
$$H=(x_2^kL_1+x_1^kL_2)_{|X},$$
where $\deg L_i\geq 1$.
Let $y\equiv(0:0:0:0:1)$, we can assume without loss of generality
that   $y\not\in \Sing(X)$, $T_yX=(x_1+x_2)$ and $L_1(y)\neq 0$, $L_2(y)\neq0$. 
The equation of $X$ is of the form 
$$(x_0x_4+x_3^2)^2+x_1F_1+x_2F_2=0,$$
to express  $X$ at the point $y$, in a better way, we can rewrite it as follows
$$x_4^3(x_1+x_2)+x_4^2(x_0^2+x_1R_1+x_2R_2)+x_4C+D=0.$$
Let $F=(x_2^kL_1+x_1^kL_2=0)$ I claim that due to the  monomial $x_0^2x_4^2$
$$\mult_y F_{|X}\leq k+1\leq \deg F.$$
Indeed let $\nu:Y\to\P^4$ be the blow up of the point $y$. 
Let $y_i$ be the coordinates in the exceptional divisor $E_0$ of equation
$(x_3=0)$ in the affine piece $y_3\neq 0$. Then
$$F_Y=(y_2^kL_1^{\prime}+y_1^kL_2^{\prime}=0){\rm,\  and}\ \ 
X_Y=(y_1+y_2+(y_0^2+y_1R_1^{\prime}+y_2R_2^{\prime})x_3+G^{\prime}x_3^2),$$
and
$$\mult_y F=k.$$
Let $\mu:W\to Y$ be the blow up of $G=X_{Y|E_0}$, with equations $x_3=y_1+y_2=0$.
Let $t$ be the coordinate in the exceptional divisor $E_1$
of equation $(x_3=0)$.
The polynomial $L_i$ does not vanish at $y$ then $F_{Y|E_0}=(\alpha y_1^k+\beta y_2^k)$, for some non zero
numbers $\alpha$ and $\beta$. 
Therefore $(y_1+y_2)^2$ does not divide $(\alpha y_1^k+\beta y_2^k)$,
and
$$\mult_{G} F_Y\leq 1.$$ 
If $\mult_{G} F_Y= 1$ then
$$F_W=(tA_0+y_1A_1+y_2A_2+x_3 B)\ \ {\rm and}\ \ X_W=(t+y_0^2+M+Nx_3),$$ 
for non zero polynomials $A_i$, and $B$, and due to the presence of the monomial $y_0^2$ the divisor
$F_{W|E_1}$ does not contain $X_{W|E_1}$ and consequently $\mult_x F_{|X}\leq k+1$.
This inequality concludes the proof.
\end{proof}
\begin{rem} I proved that a double conic is never the center
of maximal singularities on any terminal $\Q$-factorial quartic.
This relax the assumptions in \cite[Theorem 1.1]{CMp}.
\end{rem}

It is still left to adapt the proof to arbitrarily fields of characteristic 0.
\begin{proof}[Proof of Theorem \ref{th:k}] 
 Let again $\Gamma$ be a center of maximal singularities for the linear system $\H\subset|\O(n)|$.
If $\Gamma$ is defined over $k$ then all the proof works exactly as in the algebraically closed
field case. The only observations I want to add are the following.
When $\Gamma$ is a twisted cubic then $\Pi\iso\P^3\supset \Gamma$ is defined over $k$. Moreover
$H=\Pi_{|X}$ has to be smooth on the generic point of $\Gamma$, as in the proof of Lemma
\ref{le:bound}, and hence irreducible, by $\Q$-factoriality.
When $\Gamma$ is a plane curve of
degree greater than 1, the plane $\Pi\supset\Gamma$ is defined over $k$, and
$\Pi\cap X$ is a curve.

Assume that $\Gamma$ is not defined over $k$, and let $r=\deg[k(\Gamma):k]$. 
If $\Gamma=P$ is a point then  $4n^2=\H^2\cdot\O(1)\geq r(\mult_P\H)^2$.
If $P$ is smooth then $\mult_P\H^2>4n^2$, \cite[Theorem 3.1]{Co2}, 
while for singular $P$, $\mult_P\H>n$, \cite[Theorem 3.10]{Co2}, and consequently
$\mult_P\H^2>2n^2$, the exceptional divisor is a quadric.
This proves that when $r\geq 2$ no point can be a center of maximal singularities.

If $\Gamma$ is a curve then again by numerical reasons 
$$4n=\H\cdot\O(1)^2\geq r\deg\Gamma\mult_{\Gamma}\H>2n\deg\Gamma,$$
 so that
$\Gamma$ is a line and $r\leq 3$.
Let $\Gamma_i$ the conjugate lines over $k$. First observe that $\Gamma\cap\Gamma_i\neq \emptyset$.
Indeed they are both centers of maximal singularities on $\overline{k}$
and we can untwist $\Gamma$ over $\overline{k}$.
If $\Gamma_i$ is disjoint from $\Gamma$ the untwist is an isomorphism on
the generic point of $\Gamma_i$. This is very clear from our description in terms of Sarkisov
links. Then
its strict transform is a curve, say $\Gamma_i^{\prime}$, of degree $g>1$. 
Let $\H'$ be the untwist of $\H$, then 
by Lemma \ref{le:untw}, $\H'\in|\O(n')|$, for some $n'<n$.  But then $\mult_{\Gamma_i^{\prime}}\H'>n$
and this is not allowed by the proof of Theorem \ref{th:cA1}.
 
To conclude we have to study conjugate lines intersecting in a point.
Assume that $r=2$, denote $\Pi \subset \P^4$ the plane spanned by $\Gamma$ and $\Gamma_1$.
Let $S$,
$S'$ be general members of the linear system $|\I_{\Pi,X}(1)|$.
Observe that $\Pi$ is defined over $k$, therefore $\Pi\cap X=(\Gamma+\Gamma_1)+
\Delta$ is a curve.
By the proof of Theorem \ref{th:center}, since $X$ has only ordinary double points,
 all singular points of $S$ are of type
$0\in(xy+z^{t+1}=0)$, with $t\leq 2$. 

Assume that $\Pi_{|X}\neq 2(\Gamma+\Gamma_1)$  then
following the same arguments of page \pageref{pag:plane} we have to prove that for any irreducible
curve $C\subset\Delta$
$$(\Gamma+\Gamma_1)\cdot C\geq \deg C$$

Fix a point  $p\in C\cap\Gamma$.
Since both $C$ and $\Gamma$ are curves contained in $\Pi$ and $p\in S\sim(0\in (xy+z^t=0))$, with $t\leq 3$,
then 
$$(C\cdot\Gamma)_x\geq \frac12.$$
Similarly for $\Gamma_1$, so that
$$C\cdot(\Gamma+\Gamma_1)\geq 2\frac12 \deg C\geq\deg C.$$ 

If $\Pi\cap X=2(\Gamma+\Gamma_1)$ then, up to a projectivity,
 we can write the equation of
$X/\overline{k}$ as
$$x_0^2x_4^2+x_1F_1+x_2F_2=0.$$
Then we derive a contradiction as in the double conic case. Keep in mind that the crucial
point was the presence of the monomial $x_0^2x_4^2$. 

Note that the two lines are centers 
of maximal singularities on $\overline{k}$. 
Here we proved that they are not centers of maximal singularities
 with the same associated linear system.
The case $r=3$ is similar. If all lines stays on the same $k$-plane
I conclude as above. If they span a $\P^3$ say $\Pi$, then $\Pi$ is defined over $k$.
Moreover $H=\Pi_{|X}$ 
has to be smooth on the generic point of the lines, as in the proof of Lemma
\ref{le:bound}. Therefore the plane spanned by each pair of lines is not contained in $X$ and I conclude
as before.
\end{proof}

\end{document}